\magnification\magstep1
\openup1\jot

\catcode`\À=\active \defÀ{\`A}    \catcode`\à=\active \defà{\`a} 
\catcode`\Â=\active \defÂ{\^A}    \catcode`\â=\active \defâ{\^a} 
\catcode`\Æ=\active \defÆ{\AE}    \catcode`\æ=\active \defæ{\ae}
\catcode`\Ç=\active \defÇ{\c C}   \catcode`\ç=\active \defç{\c c}
\catcode`\È=\active \defÈ{\`E}    \catcode`\è=\active \defè{\`e} 
\catcode`\É=\active \defÉ{\'E}    \catcode`\é=\active \defé{\'e} 
\catcode`\Ê=\active \defÊ{\^E}    \catcode`\ê=\active \defê{\^e} 
\catcode`\Ë=\active \defË{\"E}    \catcode`\ë=\active \defë{\"e} 
\catcode`\Î=\active \defÎ{\^I}    \catcode`\î=\active \defî{\^\i}
\catcode`\Ï=\active \defÏ{\"I}    \catcode`\ï=\active \defï{\"\i}
\catcode`\Ô=\active \defÔ{\^O}    \catcode`\ô=\active \defô{\^o} 
\catcode`\Ù=\active \defÙ{\`U}    \catcode`\ù=\active \defù{\`u} 
\catcode`\Û=\active \defÛ{\^U}    \catcode`\û=\active \defû{\^u} 
\catcode`\Ü=\active \defÜ{\"U}    \catcode`\ü=\active \defü{\"u} 

\catcode`\ =\active \def { }
\def\bquote{``\thinspace\thinspace\bgroup\let = }
\def\equote{\egroup\thinspace\thinspace''}

\hsize=11.25cm    
\vsize=18cm       
\parindent=12pt   \parskip=5pt     

\hoffset=.5cm   
\voffset=.8cm   

\pretolerance=500 \tolerance=1000  \brokenpenalty=5000

\catcode`\@=11

\font\eightrm=cmr8         \font\eighti=cmmi8
\font\eightsy=cmsy8        \font\eightbf=cmbx8
\font\eighttt=cmtt8        \font\eightit=cmti8
\font\eightsl=cmsl8        \font\sixrm=cmr6
\font\sixi=cmmi6           \font\sixsy=cmsy6
\font\sixbf=cmbx6

\font\tengoth=eufm10 
\font\eightgoth=eufm8  
\font\sevengoth=eufm7      
\font\sixgoth=eufm6        \font\fivegoth=eufm5

\skewchar\eighti='177 \skewchar\sixi='177
\skewchar\eightsy='60 \skewchar\sixsy='60

\newfam\gothfam           \newfam\bboardfam

\def\tenpoint{
  \textfont0=\tenrm \scriptfont0=\sevenrm \scriptscriptfont0=\fiverm
  \def\rm{\fam\z@\tenrm}
  \textfont1=\teni  \scriptfont1=\seveni  \scriptscriptfont1=\fivei
  \def\oldstyle{\fam\@ne\teni}\let\old=\oldstyle
  \textfont2=\tensy \scriptfont2=\sevensy \scriptscriptfont2=\fivesy
  \textfont\gothfam=\tengoth \scriptfont\gothfam=\sevengoth
  \scriptscriptfont\gothfam=\fivegoth
  \def\goth{\fam\gothfam\tengoth}
  
  \textfont\itfam=\tenit
  \def\it{\fam\itfam\tenit}
  \textfont\slfam=\tensl
  \def\sl{\fam\slfam\tensl}
  \textfont\bffam=\tenbf \scriptfont\bffam=\sevenbf
  \scriptscriptfont\bffam=\fivebf
  \def\bf{\fam\bffam\tenbf}
  \textfont\ttfam=\tentt
  \def\tt{\fam\ttfam\tentt}
  \abovedisplayskip=12pt plus 3pt minus 9pt
  \belowdisplayskip=\abovedisplayskip
  \abovedisplayshortskip=0pt plus 3pt
  \belowdisplayshortskip=4pt plus 3pt 
  \smallskipamount=3pt plus 1pt minus 1pt
  \medskipamount=6pt plus 2pt minus 2pt
  \bigskipamount=12pt plus 4pt minus 4pt
  \normalbaselineskip=12pt
  \setbox\strutbox=\hbox{\vrule height8.5pt depth3.5pt width0pt}
  \let\bigf@nt=\tenrm       \let\smallf@nt=\sevenrm
  \normalbaselines\rm}

\def\eightpoint{
  \textfont0=\eightrm \scriptfont0=\sixrm \scriptscriptfont0=\fiverm
  \def\rm{\fam\z@\eightrm}
  \textfont1=\eighti  \scriptfont1=\sixi  \scriptscriptfont1=\fivei
  \def\oldstyle{\fam\@ne\eighti}\let\old=\oldstyle
  \textfont2=\eightsy \scriptfont2=\sixsy \scriptscriptfont2=\fivesy
  \textfont\gothfam=\eightgoth \scriptfont\gothfam=\sixgoth
  \scriptscriptfont\gothfam=\fivegoth
  \def\goth{\fam\gothfam\eightgoth}
  
  \textfont\itfam=\eightit
  \def\it{\fam\itfam\eightit}
  \textfont\slfam=\eightsl
  \def\sl{\fam\slfam\eightsl}
  \textfont\bffam=\eightbf \scriptfont\bffam=\sixbf
  \scriptscriptfont\bffam=\fivebf
  \def\bf{\fam\bffam\eightbf}
  \textfont\ttfam=\eighttt
  \def\tt{\fam\ttfam\eighttt}
  \abovedisplayskip=9pt plus 3pt minus 9pt
  \belowdisplayskip=\abovedisplayskip
  \abovedisplayshortskip=0pt plus 3pt
  \belowdisplayshortskip=3pt plus 3pt 
  \smallskipamount=2pt plus 1pt minus 1pt
  \medskipamount=4pt plus 2pt minus 1pt
  \bigskipamount=9pt plus 3pt minus 3pt
  \normalbaselineskip=9pt
  \setbox\strutbox=\hbox{\vrule height7pt depth2pt width0pt}
  \let\bigf@nt=\eightrm     \let\smallf@nt=\sixrm
  \normalbaselines\rm}

\tenpoint

\def\pc#1{\bigf@nt#1\smallf@nt}         \def\pd#1 {{\pc#1} }

\catcode`\;=\active
\def;{\relax\ifhmode\ifdim\lastskip>\z@\unskip\fi
\kern\fontdimen2  -1.2 \fontdimen3 \string;}

\catcode`\:=\active
\def:{\relax\ifhmode\ifdim\lastskip>\z@\unskip\fi\penalty\@M\ \fi\string:}

\catcode`\!=\active
\def!{\relax\ifhmode\ifdim\lastskip>\z@
\unskip\fi\kern\fontdimen2  -1.1 \fontdimen3 \string!}

\catcode`\?=\active
\def?{\relax\ifhmode\ifdim\lastskip>\z@
\unskip\fi\kern\fontdimen2  -1.1 \fontdimen3 \string?}

\catcode`\«=\active 
\def«{\raise.4ex\hbox{%
 $\scriptscriptstyle\langle\!\langle$}}

\catcode`\»=\active 
\def»{\raise.4ex\hbox{%
 $\scriptscriptstyle\rangle\!\rangle$}}

\frenchspacing

\def\raggedbottom{\topskip 10pt plus 36pt\r@ggedbottomtrue}

\def\pointir{\unskip . --- \ignorespaces}

\def\Medbreak{\vskip-\lastskip\medbreak}

\long\def\th#1 #2\enonce#3\endth{
   \Medbreak\noindent
   {\pc#1} {#2\unskip}\pointir{\it #3}\smallskip}

\def\proof{\vskip-\lastskip\smallskip\noindent
 {\it Proof} : }

\def\decale#1{\smallbreak\hskip 28pt\llap{#1}\kern 5pt}
\def\decaledecale#1{\smallbreak\hskip 34pt\llap{#1}\kern 5pt}
\def\puce{\smallbreak\hskip 6pt{$\scriptstyle\bullet$}\kern 5pt}

\def\eqalign#1{\null\,\vcenter{\openup\jot\m@th\ialign{
\strut\hfil$\displaystyle{##}$&$\displaystyle{{}##}$\hfil
&&\quad\strut\hfil$\displaystyle{##}$&$\displaystyle{{}##}$\hfil
\crcr#1\crcr}}\,}

\catcode`\@=12

\showboxbreadth=-1  \showboxdepth=-1

\mathcode`A="7041 \mathcode`B="7042 \mathcode`C="7043 \mathcode`D="7044
\mathcode`E="7045 \mathcode`F="7046 \mathcode`G="7047 \mathcode`H="7048
\mathcode`I="7049 \mathcode`J="704A \mathcode`K="704B \mathcode`L="704C
\mathcode`M="704D \mathcode`N="704E \mathcode`O="704F \mathcode`P="7050
\mathcode`Q="7051 \mathcode`R="7052 \mathcode`S="7053 \mathcode`T="7054
\mathcode`U="7055 \mathcode`V="7056 \mathcode`W="7057 \mathcode`X="7058
\mathcode`Y="7059 \mathcode`Z="705A

\font\ss=cmss10 

\def\qp{{\bf Q}_p}
\def\qmodz{{\bf Q}/{\bf Z}}
\def\azero#1{\hbox{{\ss A}}_0(#1)}
\def\azeroo#1{\azero{#1}_0}
\def\zzero#1{\hbox{{\ss Z}}_0(#1)}
\def\zzeroo#1{\zzero{#1}_0}
\def\br#1{\hbox{{\ss Br}}(#1)}
\def\brzero#1{\br{#1}_0}
\def\acheun{\hbox{\ss H}^1}

\def\Abar{\overline A}
\def\Kbar{{\overline K}}
\def\Pbar{{\overline P}}
\def\Lbar{{\overline L}}
\def\Xbar{\overline X}
\def\Ubar{\overline U}
\def\Xprim{{X'}}

\def\Abaretoile{{\overline A}{}^\times}

\def\Letoile{L^\times}

\def\xibar{\overline{\xi}}
\def\etabar{\overline{\eta}}

\def\Ketoile{K^{\times}}
\def\Letoile{L^{\times}}

\def\Kbaretoile{{\Kbar{}^{\times}}}
\def\Lbaretoile{{\Lbar{}^{\times}}}
\def\Kprim{{K'}}

\def\Gal{\mathop{\hbox{\ss Gal}}\nolimits}
\def\Pic{\mathop{\hbox{\ss Pic}}\nolimits}

\def\Hom{\mathop{\hbox{\ss Hom}}\nolimits}
\def\Ext{\mathop{\hbox{\ss Ext}}\nolimits}
\def\extun{\Ext^1}
\def\extensionde{\,|\,}
\def\End{\mathop{\rm End}\nolimits}

\def\Spec{\mathop{\hbox{\ss Spec}}\nolimits}

\def\P{\mathord{\bf P}}
\def\Z{\mathord{\bf Z}}
\def\Q{\mathord{\bf Q}}

\def\hfl#1#2#3{\smash{\mathop{\hbox to#3{\rightarrowfill}}\limits
^{\textstyle#1}_{\textstyle#2}}}
\def\gfl#1#2#3{\smash{\mathop{\hbox to#3{\leftarrowfill}}\limits
^{\textstyle#1}_{\textstyle#2}}}

\def\pafl#1#2#3{\llap{$\textstyle #1$}
\left\Vert\vbox to#3{}\right.\rlap{$\textstyle #2$}}

\def\qed{\raise -2pt\hbox{\vrule\vbox to 10pt{\hrule width 4pt
                 \vfill\hrule}\vrule}}

\def\phi{\varphi}

\def\Id{\mathop{\rm Id}\nolimits}
\def\Ind{\mathop{\hbox{\ss Ind}}\nolimits}
\def\droite#1{\,\hfl{#1}{}{8mm}\,}
\def\versbas#1{\vfl{#1}{}{5mm}}
\def\vershaut#1{\ufl{#1}{}{5mm}}
\def\cores{\mathop{\hbox{\ss Cor}}\nolimits}

\def\diagram#1{\def\normalbaselines{\baselineskip=0pt\lineskip=5pt}
\matrix{#1}}
\def\vfl#1#2#3{\llap{$\textstyle #1$}
\left\downarrow\vbox to#3{}\right.\rlap{$\textstyle #2$}}

\def\pafl#1#2#3{\llap{$\textstyle #1$}
\left\Vert\vbox to#3{}\right.\rlap{$\textstyle #2$}}

\def\ufl#1#2#3{\llap{$\textstyle #1$}
\left\uparrow\vbox to#3{}\right.\rlap{$\textstyle #2$}}

\newcount\refno 

\long\def\ref#1:#2<#3>{                                        
\global\advance\refno by1\par\noindent                              
\llap{[{\bf\number\refno}]\ }{#1} \pointir{\it #2} #3\goodbreak }

\def\citer#1(#2){[{\bf\number#1}\if#2\empty\relax\else,\ #2\fi]}

\def\fleche{\rightarrow}

\newcount\numerodesection
\def\section#1{\bigbreak
 {\bf\number\numerodesection.\ \ #1}\nobreak\medskip
 \advance\numerodesection by1}

\newcount\numeroderemarque
\def\remarque{\advance\numeroderemarque by1\smallbreak
{\it Remarque\/}\ \number\numeroderemarque~:}

\def\remark#1{\smallbreak\noindent
{\it Remark\/} {#1\unskip}\pointir}

\newcount\formuleno
\def\numeroter{\global\advance\formuleno by1
 \leqno{(\oldstyle\number\formuleno)}}
\def\formule#1{$(\oldstyle\number#1)$}

\def\somme{\mathop{\smash{\raise 2pt\hbox{$\sum$}}}\limits}

\def\Arond{{\cal A}}

\def\Orond{{\cal O}}

\def\zero{\{0\}}
\def\one{\{1\}}
\def\long{\mathop{\rm long}\nolimits}

\newbox\bibbox
\setbox\bibbox\vbox{\bigbreak
\centerline{{\pc BIBLIOGRAPHICAL REFERENCES}}

\ref{\pc ARTIN} (M.), {\pc GROTHENDIECK} (A.) and {\pc VERDIER} (J.-L.):
Th\'eorie des topos et cohomologie \'etale des sch\'emas,
<t.~3, S\'eminaire de G\'eom\'etrie Alg\'ebrique du Bois-Marie
1963--1964 (SGA 4), Lecture Notes in Mathematics, vol.\ 305,
Springer-Verlag, Berlin, 1973.>
\newcount\sgaiv \global\sgaiv=\refno

\ref{\pc BLOCH} (S.):
On the Chow groups of certain rational surfaces,
<Ann.\ Sci.\ \'Ecole Norm.\ Sup.\ (4) {\bf 14} (1981) 1, 41--59.>
\newcount\bloch \global\bloch=\refno

\ref {\pc COLLIOT}-{\pc THÉLÈNE} (J.-L.):
Hilbert's theorem $90$ for $K_2$, with applications to Chow groups of
rational surfaces, 
<Invent.\ math.\ {\bf 71} (1983), 1--20.>
\newcount\hilbertxc \global\hilbertxc=\refno

\ref{\pc COLLIOT}-{\pc TH{\'E}L{\`E}NE} (J.-L.):
Lettre à C.~S.~Dalawat,
<24 octobre 1998.>
\newcount\lettrecolliot \global\lettrecolliot=\refno

\ref{\pc COLLIOT}-{\pc THÉLÈNE} (J.-L.), {\pc CORAY} (D.) and 
{\pc SANSUC} (J.-J.): 
Descente et principe de {H}asse pour certaines vari\'et\'es
rationnelles,
<J.\ Reine Angew.\ Math.\ {\bf 320} (1980), 150--191.>
\newcount\ctcoraysansuc \global\ctcoraysansuc=\refno

\ref{\pc COLLIOT}-{\pc THÉLÈNE} (J.-L.) and {\pc SAITO} (S.):
Z\'ero-cycles sur les vari\'et\'es $p$-adiques et groupe de
Brauer,
<Internat.\ Math.\ Res.\ Notices (1996) 4, 151--160.>
\newcount\ctsaito \global\ctsaito=\refno

\ref{\pc COLLIOT}-{\pc THÉLÈNE} (J.-L.) and {\pc SANSUC} (J.-J.): 
On the Chow groups of certain rational
surfaces : a sequel to a paper of S.~Bloch, 
<Duke math.\ jour.\ {\bf 48} (1981), 421--427.>
\newcount\ctssequel \global\ctssequel=\refno

\ref{\pc CORAY} (D. F.) and {\pc TSFASMAN} (M. A.):
Arithmetic on singular del {P}ezzo surfaces,
<Proc.\ London Math.\ Soc.\ (3) {\bf 57} (1), 1988, p.~25--87.>
\newcount\coraytsfasman  \global\coraytsfasman=\refno

\ref{\pc DALAWAT} (C. S.):
Le groupe de Chow d'une surface de Châtelet sur un corps local,
<Indag.\ mathem.\ N.S. {\bf 11} (2) (2000), 173--185. {\tt math.AG/0302156}.>
\newcount\chatelet  \global\chatelet=\refno

\ref{\pc DALAWAT} (C. S.):
Le groupe de Chow d'une surface rationnelle sur un corps local,
<{\tt math.AG/0302157}.>
\newcount\chow  \global\chow=\refno

\ref{\pc FULTON} (W.):
Intersection Theory,
<Springer, Berlin, 1984.>
\newcount\fulton  \global\fulton=\refno


\ref{\pc GROTHENDIECK} (A.):
Le groupe de Brauer\/ {\rm I, II, III, 46--188} 
<dans Dix exposés sur la cohomologie des schémas, 
North-Holland, Amsterdam, 1968.>
\newcount\grothendieck  \global\grothendieck=\refno

\ref{\pc LICHTENBAUM} (S.):
Duality theorems for curves over $p$-adic fields,
<Invent.\ math.\ {\bf 7} (1969), 120--136.>
\newcount\lichtenbaum  \global\lichtenbaum=\refno

\ref{\pc MANIN} (Y. I.):
Le groupe de Brauer-Grothendieck en g\'eom\'etrie diophantienne,
<Actes du Congr\`es International des Math\'ematiciens (Nice, 1970),
t.~1, 401--411, Gauthier-Villars, Paris, 1971.>
\newcount\nice  \global\nice=\refno

\ref{\pc MANIN} (Y. I.):
Cubic forms, 
<$2^{\rm nd}$ edition, North-Holland, Amsterdam, 1986.>
\newcount\cubic  \global\cubic=\refno

\ref{\pc MILNE} (J.):
The Brauer group of a rational surface,
<Invent.\ math.\ {\bf 11} (1970), 304--307.>
\newcount\milne  \global\milne=\refno

\ref{\pc SANSUC} (J.-J.):
{\`A} propos d'une conjecture arithm{\'e}tique sur le groupe de
{C}how d'une surface rationnelle,
<S{\'e}minaire de Th{\'e}orie des nombres de Bordeaux, Expos{\'e} 33, 1982.>
\newcount\sansuc  \global\sansuc=\refno

\ref {\pc SERRE} (J.-P.):
Galois Cohomology,
<Springer, Berlin, 1997.>
\newcount\serrecg  \global\serrecg=\refno

} 


\centerline{\bf Restriction, corestriction}     
\smallskip
\centerline{\bf and }
\smallskip
\centerline{\bf the characteristic homomorphism}
\smallskip
\centerline{\bf for}
\smallskip
\centerline{\bf $0$-cycles of degree~$0$}

\vskip4mm

\centerline{Chandan Singh {\pc DALAWAT}}

\vskip2cm

\section{Introduction}
The main purpose of this Note is to verify (see \S\thinspace5) that
the characteristic homomorphism \citer\ctssequel(p.~423) of
J.-L. Colliot-Th{\'e}l{\`e}ne and J.-J. Sansuc is compatible with
restriction and corestriction maps on the Chow groups of $0$-cycles of
degree~$0$.  We also derive some consequences of this compatibility
(see \S\S\thinspace3, 6) ; as a matter of fact, it is these
applications which triggered this investigation.

A word about the genesis and the organisation of this Note is in
order.  In the middle of February~2003, for reasons which need not
concern us here, I decided to write up for publication a letter of
Colliot-Th{\'e}l{\`e}ne \citer\lettrecolliot()
(cf.~\citer\chatelet(remarque~2)) in which he had indicated how the
explicit proof of \citer\chatelet(prop.~3) --- for Ch{\^a}telet
surfaces --- was a consequence of a proposition (see prop.~1.1) valid
for all rational surfaces over a field $K$ which is a finite
extensions of $\qp$ ($p$ a prime number).  The process of writing led
me to an equivalent formulation (see remark~4.3) --- perhaps a simpler
one --- in terms of $\Ext$ groups in place of $\hbox{\ss Br}$auer
groups ; it transpired that this formulation retains its validity over an
arbitrary perfect field $K$ 
(see \S\thinspace5) ; it leads to a similar statement for rational
surfaces having a $K$-point ; the case of a number fields (where the
existence of a $K$-point is not needed) can be found in \S\thinspace6.
The first three sections (\S\S\thinspace1--3), which treat the local
case, reproduce --- with his kind permission and with no significant
change --- the contents of
\citer\lettrecolliot().  The \S\S\thinspace4--6 treat the general
case, including number fields.

{\it Je remercie Jean-Louis Colliot-Thélène pour sa
lettre}~\citer\lettrecolliot(), {\it qui est à l'origine de cette
Note, et pour ses conseils, ses commentaires et sa patience.  I thank
Joost van Hamel for his interest in this work, and for a critical
reading of the manuscript.}

\section{The local statement}

Let us take the quickest route to the statement of the main
proposition of \citer\lettrecolliot() ; some of the definitions will
be recalled in greater detail in the course of the proof in
\S\thinspace2.

For a variety $X$ over a field $K$, we shall denote by $\br{X}$
the group of equivalence classes of $\Orond_X$-algebras which are
locally isomorphic to the $\Orond_X$-algebra $\End_{\Orond_X} V$
of endomorphisms of a vector bundle $V$ over $X$, modulo those
algebras which are globally isomorphic to such an algebra of
endomorphisms (the Brauer group \citer\grothendieck())~; the fibre of
such an algebra at any $x\in X$ is a central simple $K(x)$-algebra.
When $X$ is regular of dimension $\le2$, the natural injection
$\br{X}\fleche\hbox{\ss H}^2(X,{\bf G}_m)$ is an isomorphism
\citer\grothendieck(Brauer~II, cor.~2.2).  We
write $\br{K}$ instead of $\br{\Spec K}$ and $\brzero{X}$
for the quotient $\br{X}/\br{K}$ modulo the subgroup of
\bquote constant\equote\ algebras ; this will be called the
{\it reduced\/} Brauer group of $X$.  (Note that the map
$\br{K}\fleche\br{X}$ need not be injective, even for curves
\citer\lichtenbaum()).  The group of $0$-cycles on $X$, modulo
rational equivalence, will be denoted by $\azero X$ and --- when $X$
is proper --- $\azeroo X$ will stand for the subgroup of degree-$0$
$0$-cycles (the {\it reduced\/} Chow group, see for example
\citer\fulton()).

Let $L$ be a finite extension of $K$ and put $X_L=X\times_K L$.  We
have the {\it restriction maps}
$$
f^*:\azeroo{X}\droite{}\azeroo{X_L}\,,\qquad 
f^*:\brzero{X}\droite{}\brzero{X_L}
\numeroter \newcount\resdef \global\resdef=\formuleno
$$
and the  {\it corestriction maps}
$$
f_*:\azeroo{X_L}\droite{}\azeroo{X}\,,\qquad 
f_*:\brzero{X_L}\droite{} \brzero{X}.
\advance\abovedisplayskip by5pt
\numeroter \newcount\coresdef
\global\coresdef=\formuleno
\advance\belowdisplayskip by-\bigskipamount
$$

\th PROPOSITION 1.1 (Colliot-Th{\'e}l{\`e}ne)
\enonce
Let\/ $L\extensionde K\extensionde \qp$ ($p$ prime) be finite extensions, and let\/
$X$ be a smooth, proper, absolutely connected\/ $K$-surface
potentially birational to\/ $\P_2$.

a) If the restriction~$f^*$ \formule\resdef\ on reduced Brauer groups
is surjective (resp.~$0$), then the corestriction~$f_*$
\formule\coresdef\ on reduced Chow groups is injective (resp.~$0$).

b) If the corestriction~$f_*$ \formule\coresdef\ on reduced Brauer
groups is surjective (resp.~$0$), then the restriction~$f^*$
\formule\resdef\ on reduced Chow groups is injective (resp.~$0$).
\endth

%
 
\section{Generalities and proof}

For a variety $X$ over a field $K$, we denote by $\zzero{X}$ the free
commutative group on the set of closed points of $X$ and by $\zzeroo{X}$
the subgroup of $0$-cycles of degree~$0$.  Manin (\citer\nice(),
\citer\cubic()) defined a natural bilinear pairing
$$
\langle\phantom{x},\phantom{a}\rangle\;:\;
\zzero{X}\times\br{X}\droite{}\br{K}
$$
which to a closed point $x\in X$ --- of residue field $K'=K(x)$ ---
and a class $a\in\br{X}$, associates the class
$\cores_{K'|K}(a(x))\in\br{K}$, where $a(x)$ is the class in $\br{K'}$
of the \bquote fibre\equote\ of $a$ at $x$ (if $a$ is represented by
an $\Orond_X$-algebra $\Arond$, then $a(x)\in\br{K'}$ is represented
by the fibre $\Arond(x)=\Arond_x\otimes_{\Orond_{X,x}}K'$ of $\Arond$
at $x$, which is a central simple $K'$-algebra).  When $a\in\br{X}$
comes from an element of $\br{K}$ ({\it i.e.}  when it is a
\bquote constant\equote\ algebra), one has
$\langle z,a\rangle=\deg(z)a$ ; in particular, $\langle z,a\rangle=0$
for every $z\in\zzeroo{X}$ and for every \bquote constant\equote\
class $a\in\br{X}$.

Suppose now that the $K$-variety $X$ is proper.  Then
$\langle z,a\rangle=0$ for any $0$-cycle $z$ rationally equivalent to
$0$ and for every $a\in\br{X}$ (\citer\sgaiv(Exp.~XVII,
\S\thinspace6), cf.~\citer\bloch(Appendix), \citer\ctssequel()).  We
thus get a pairing
$$
\langle\phantom{x},\phantom{a}\rangle\;:\;
\azero{X}\times\br{X}\droite{}\br{K}
$$
where $\azero{X}$ is the group of $0$-cyles on $X$, modulo rational
equivalence. Let $\azeroo{X}$ be the kernel of the degree map
$\deg:\azero{X}\fleche\Z$ and recall that $\brzero{X}$ stands
for $\br{X}/\br{K}$.  We have an induced pairing
(cf.~\citer\ctsaito(p.~153))
$$
\langle\phantom{x},\phantom{a}\rangle\;:\;
\azeroo{X}\times\brzero{X}\droite{}\br{K}.
\numeroter \newcount\thepairingk
\global\thepairingk=\formuleno
$$
Now let $f:\Spec L\fleche\Spec K$ be the structure map corresponding
to a finite separable extension $L\extensionde K$.  Put $X_L=X\times_K L$ ; we
thus have a pairing
$$
\langle\phantom{x},\phantom{a}\rangle\;:\;
\azeroo{X_L}\times\brzero{X_L}\droite{}\br{L}
\droite{f_*}\br{K}
\numeroter \newcount\thepairingl
\global\thepairingl=\formuleno
$$
in which $f_*$ denotes the corestriction map on Brauer groups of fields.

As consequences of the functoriality of the definitions, one has the
two projection formulas
$$
\eqalign{
f_*\langle z, f^*_X(a)\rangle&=
\langle f_{X*}(z), a\rangle&\qquad\hbox{for\ } 
&z\in\zzero{X_L},\ a\in\br{X},\cr
f_*\langle f^*_X(z), a\rangle &=
\langle z, f_{X*}(a)\rangle&\qquad\hbox{for\ }   
&z\in\zzero{X},\ a\in\br{X_L},
}
\numeroter \newcount\projform \global\projform=\formuleno
$$
where $f_X:X_L\fleche X$ denotes the canonical morphism and where
$f_{X*}$ serves to denote the corestriction maps for $0$-cycles as
well as for the Brauer group and, similarly, $f^*_X$ denotes the
restriction maps.  To simplify notation, we write $f_*$, $f^*$
--- as we did in \formule\resdef, \formule\coresdef\ ---
instead of $f_{X*}$, $f^*_X$.


Now suppose that $K$ is a finite extension of $\qp$ ($p$ a prime
number).  Then both $\br{K}$ and $\br{L}$ are canonically
isomorphic to $\qmodz$ by local class field theory, and the
corestriction map $f_*:\br{L}\fleche\br{K}$ is just the
identity map $\Id:\qmodz\fleche\qmodz$.

Denoting by $(\phantom{G})^\vee$ the functor
$\Hom(\phantom{G},\qmodz)$ on the category of commutative groups, the
pairings \formule\thepairingk\ and \formule\thepairingl\ furnish
the vertical arrows in the diagrams
{\def\droite#1{\kern-5pt\hfl{#1}{}{8mm}\kern-5pt}
$$\diagram{
\azeroo{X_L}&\droite{f_*}&
\azeroo{X}\phantom{,}&\qquad&
\azeroo{X}&\droite{f^*}&
\azeroo{X_L}\phantom{,}\cr
\versbas{}&&\versbas{}&\qquad&\versbas{}&&\versbas{}\cr
\brzero{X_L}^\vee&\droite{f^*{}^\vee}&
\brzero{X}^\vee,&\qquad&
\brzero{X}^\vee&\droite{f_*^\vee}&
\brzero{X_L}^\vee,\cr 
}\numeroter \newcount\diagcomm \global\diagcomm=\formuleno
$$} 
and the two projection formulas \formule\projform\ show that they
are {\it commutative}.
 
{\it Completion of the proof of prop.}~1.1.--- Now suppose that $X$ is a
rational surface, which means that it becomes birational to $\P_2$
over a suitable (finite) extension of $K$ --- it is sufficient to assume
that the vertical arrows in \formule\diagcomm\ are injective.  Then
\citer\hilbertxc({prop.~5, cf.~prop.~7{\it b\/})})
says that the vertical arrows in \formule\diagcomm\ are injective.
Now, if $f^*:\brzero{X}\fleche\brzero{X_L}$
\formule\resdef\ is surjective, then 
${f^*{}^\vee}:\brzero{X_L}^\vee\fleche\brzero{X}^\vee$ is
injective,  therefore prop.~1.1{\it a\/}) is an immediate consequence
of \citer\hilbertxc() and the commutativity of the first square
\formule\diagcomm.  Similarly, if 
$f_*:\brzero{X_L}\fleche\brzero{X}$ \formule\coresdef\ is
surjective, then
${f_*{}^\vee}:\brzero{X}^\vee\fleche\brzero{X_L}^\vee$ is
injective, and prop.~1.1{\it b\/}) follows from \citer\hilbertxc() and
the commutativity of the second square \formule\diagcomm.

\section{The local application}

To illustrate the algebraic solution provided by prop.~1.1 to the
arithmetical problem of determining the restriction and corestriction
maps on Chow groups ($0$-cycles of degree~0 modulo rational
equivalence) of rational surfaces over local fields, we shall content
ourselves with deriving the following proposition, which appears in
\citer\chatelet(prop.~3) with the superfluous hypothesis
\bquote $p\neq2$\equote : 

\th PROPOSITION 3.1
\enonce 
Let\/ $L\extensionde K\extensionde \qp$ ($p$ prime) be finite extensions, and let\/
$X$ be a smooth projective\/ $K$-surface which is\/ $K$-birational to
$$
y^2-dz^2=x(x-e_1)(x-e_2)\quad 
(e_1\neq e_2\hbox{\rm\ in\ }\Ketoile,\ 
d\in\Ketoile,\notin\Ketoile{}^2).
\numeroter
$$ 
Then the homomorphism of restriction\/ $\azeroo{X}\fleche\azeroo{X_L}$
is trivial if the degree\/ $n=[L:K]$ is even ; it is an isomorphism
if $n$ is odd.
\endth
\newcount\chatsurf \global\chatsurf=\formuleno

The proof will occupy the rest of this \S.  An example of such an $X$
is provided by the surface \citer\chatelet({(2)}) ; it is
$K(\sqrt d)$-birational to $\P_2$.

If $L$ is not linearly disjoint from $\Kprim=K(\sqrt d)$ --- $n$ is
then even ---, the surface $X_L$ is $L$-birational to $\P_2$ ; we then
have $\azeroo{X_L}=\zero$, and there is nothing to prove.  Assume
henceforth that $L$ does not contain $\sqrt d$.

We shall show that the restriction ~$f^*$ \formule\resdef\ on reduced
Brauer groups is always an isomorphism ; the corestriction ~$f_*$
\formule\coresdef\ on reduced Brauer groups is an isomorphism if $n$
is odd and $0$ if $n$ is even.  An application of prop.~1.1 will
conclude the argument.

Let $\Lbar$ be an algebraic closure of $L$ and put
$\Xbar=X\times_K\Lbar$, $\Gamma_K=\Gal(\Lbar|K)$,
$\Gamma_L=\Gal(\Lbar|L)$~; the Picard group $\Pbar=\Pic\Xbar$ --- a
free $\Z$-module of finite rank --- carries a continuous action of
$\Gamma_K$ and of its subgroup (of finite index) $\Gamma_L$, so one
has restriction and corestriction maps
{\advance\abovedisplayskip by-,5\baselineskip
\advance\belowdisplayskip by-1,5\baselineskip
$$
f^*:\acheun(K,\Pbar)\fleche\acheun(L,\Pbar),
\qquad 
f_*:\acheun(L,\Pbar)\fleche\acheun(K,\Pbar)
\advance\abovedisplayskip by5pt
$$}
{\advance\abovedisplayskip by-\baselineskip
\advance\belowdisplayskip by-2\baselineskip
\th LEMMA 3.2
\enonce
One has the following commutative diagrams in which the vertical
arrows are isomorphisms
{\def\droite#1{\kern-5pt\hfl{#1}{}{8mm}\kern-5pt}
$$\diagram{
\brzero{X}&\droite{f^*}&
\brzero{X_L}&\quad&
\brzero{X_L}&\droite{f_*}&
\brzero{X}\phantom{.}\cr
\versbas{}&&\versbas{}&\quad&\versbas{}&&\versbas{}\cr
\acheun(K,\Pbar)&\droite{f^*}&
\acheun(L,\Pbar),&\quad&
\acheun(L,\Pbar)&\droite{f_*}&
\acheun(K,\Pbar).\cr
}\numeroter
$$}
\endth
\newcount\diagbracheun
\global\diagbracheun=\formuleno
}
\proof  That the map $\brzero{X}\fleche\acheun(K,\Pbar)$ in
\formule\diagbracheun\ is an isomorphism is just \citer\nice(lemme~3)
--- applicable because $\hbox{\ss H}^3(K,\Kbaretoile)=\zero$ --- since
$\br{\Xbar}=\zero$ (cf.~\citer\milne(p.~305)) and hence
$\br{X}\fleche\br{\Xbar}$ is the zero map.  Similarly,
$\brzero{X_L}\fleche\acheun(L,\Pbar)$ is an isomorphism.  The
commutativity of the two squares is a consequence of fuctoriality.

%
%

{
\def\cohomok{(\Ketoile\!/\Ketoile{}^2)^2}
\def\cohomol{(\Letoile\!/\Letoile{}^2)^2}
\def\ldroite#1{\kern-20pt\hfl{#1}{}{12mm}\kern-20pt}
\def\droite#1{\kern-5pt\hfl{#1}{}{8mm}\kern-5pt}
\def\\{\phantom{^2}}
\advance\abovedisplayskip by-\baselineskip
\advance\belowdisplayskip by-1,5\baselineskip
\th LEMMA 3.3
\enonce
One has the following commutative diagrams in which the vertical
arrows are injections and the lower horizontal ones are induced by the
inclusion\/
$\Ketoile\fleche\Letoile$, resp.\ by the norm map\/
$\Letoile\fleche\Ketoile$ :
$$\diagram{
\acheun(K,\Pbar)\\&\droite{f^*}&
\acheun(L,\Pbar)\\\phantom{,}&\ 
\acheun(L,\Pbar)\\&\droite{f_*}&
\acheun(K,\Pbar)\\\phantom{.}\cr
\versbas{}&&
\versbas{}&\ 
\versbas{}&&
\versbas{}\cr
\cohomok&\droite{f^*}&
\cohomol,&\ 
\cohomol&\droite{f_*}&
\cohomok.\cr
}\numeroter
$$
\endth
\newcount\diagacheuncoho
\global\diagacheuncoho=\formuleno
}
\proof One needs an explicit description of the structure of the
discrete $\Gamma_K$-module $\Pbar$, see for example \citer\sansuc(),
cf.~\citer\ctcoraysansuc(prop.~5.1).  Recalling that
$\Kprim=K(\sqrt d)$, put $G=\Gal(\Kprim|K)$ and
$\Xprim=X\times_K\Kprim$ ; the map $\Pic\Xprim\fleche\Pbar=\Pic\Xbar$
is an isomorphism and hence $\acheun(K,\Pbar)=\acheun(G,\Pic\Xprim)$.
Now the $G$-module $\Pic\Xprim$ is isomorphic --- up to addition of
permutation modules \citer\ctcoraysansuc(p.~177) --- to
$(\Z[G]/\Z)^2$ : the map $\Z\fleche\Z[G]$ sends $1$ to $1+\sigma$,
where $\sigma$ is the generator of $G$.  Therefore
$\acheun(K,\Pbar)=(\Z/2\Z)^2$, and, with this identification, the
extreme vertical arrows in
\formule\diagacheuncoho\ send $1\in\Z/2\Z$ to the class of
$d\in\Ketoile$ modulo $\Ketoile{}^2$.

As $L$ is by assumption linearly disjoint from $\Kprim$, a similar
result is valid for the $\Gamma_L$-module $\Pbar$, for the group
$\acheun(L,\Pbar)$ and for the middle vertical arrows in
\formule\diagacheuncoho.

{\it Completion of the proof of prop.}~3.1.--- Notice first that the
restriction $f^*:\acheun(K,\Pbar)\fleche\acheun(L,\Pbar)$ is always an
isomorphism (lemma~3.3), since $d\in\Ketoile$ \formule\chatsurf\ does
not become a square in $\Letoile$ by hypothesis.  So the corestriction
$f_*:\azeroo{X_L}\fleche\azeroo{X}$ is injective (lemma~3.2,
prop.~1.1{\it a})).

As for the corestriction
$f_*:\acheun(L,\Pbar)\fleche\acheun(K,\Pbar)$, it is induced (lemma~3.3)
by the norm map $\Letoile\fleche\Ketoile$, which sends $d$
\formule\chatsurf\ to $d^n$ (where $n=[L:K]$).  Therefore the restriction
$f^*:\azeroo{X}\fleche\azeroo{X_L}$ is injective if $n$ is odd and $0$
if $n$ is even (lemma~3.2, prop.~1.1{\it b})).  The proof is complete
if $n$ is even ; if $n$ is odd, restriction
$f^*:\azeroo{X}\fleche\azeroo{X_L}$ and corestriction
$f_*:\azeroo{X_L}\fleche\azeroo{X}$ are injective, and, as the two
groups are finite, $f^*$ and $f_*$ are both isomorphisms. 


\th COROLLARY 3.4
\enonce 
With the notation of prop.~$3.1$, the corestriction map
$f_*:\azeroo{X_L}\fleche\azeroo{X}$ is injective for all\/ $n$ ; it is 
bijective if\/ $n$ is odd.
\endth
\proof If $L$ contains $\sqrt d$, this is trivially true, since
$\azeroo{X_L}=\zero$.  Otherwise, we have seen --- in the course of
the proof of prop.~3.1 --- that this is so.

\remark{3.5} The conclusions of prop.~3.1 and of cor.~3.4
hold --- with a similar proof --- for $X$ a smooth proper surface
$K$-birational to
$$
y^2-dz^2=x(x-e_1)(x-e_2)\cdots(x-e_r)\quad(r\ge2)
\numeroter\newcount\fibcon \global\fibcon=\formuleno
$$ 
where the $e_i$ are distinct elements of\/ $\Ketoile$ and where
$d\in\Ketoile$ is not a square.

\remark{3.6} The groups $\azeroo{X}$ and $\azeroo{X_L}$ 
have been calculated for the surfaces \formule\chatsurf ; see
\citer\coraytsfasman(prop.~4.7), cf.~\citer\chatelet(prop.~1) when the
extension $\Kprim\extensionde K$ (resp. $L\Kprim\extensionde L$) is
unramified ; see \citer\chatelet(prop.~2) when these extensions are
ramified.  The values are expressed in
\citer\chatelet(prop.~4, 5) in terms of the type of possible bad
reduction of the surfaces in question.  It seems that S.~Bloch has
expressed the hope that, more generally, for all rational surfaces
over a local field, the value of the Chow group of $0$-cycles of
degree~$0$ depends only on the type of bad reduction of that surface ;
when there is good reduction, see \citer\hilbertxc({th.~A(iii)}),
cf. \citer\chow(th.~3).  But even the case of surfaces
\formule\fibcon\ has not yet been treated, as far as I know.

\goodbreak
\section{An equivalent formulation}

Prop.~1.1 can be reformulated in terms of the characteristic
homomorphism \citer\ctssequel(p.~423) of Colliot-Th{\'e}l{\`e}ne and
Sansuc ; we shall show in the next~\S\ that this homomorphism is
compatible with restriction and corestriction for more general $K$ and
$X$.

Let $S$ be the $K$-torus whose $\Gamma_K$-module of\/ $\Lbar$-rational
points is $S(\Lbar)=\Hom_{\Z}(\Pbar,\Lbaretoile)$.  Then one has
$\extun_K(\Pbar,\Lbaretoile)=\acheun(K,S(\Lbar))$, and the cup product
furnishes a duality of finite commutative groups
$$
\acheun(K,\Pbar)\times \acheun(K,S(\Lbar))\droite{}\qmodz
$$
\citer\serrecg({p.~102, th.~6(b)}).  We have thus a sequence of
isomorphisms (cf.~\S\thinspace3)
$$
\brzero{X}^\vee
=\acheun(K,\Pbar)^\vee
=\acheun(K,S(\Lbar))
=\extun_K(\Pbar,\Lbaretoile).
$$
The same holds at the level of $L$ instead of $K$.  As these
identifications are compatible with restriction and corestriction,
they 
yield two commuatative squares 
--- which we coalesce for typographical reasons ---
{\def\droite#1{\kern-5pt\hfl{#1}{}{8mm}\kern-5pt}
\def\ldroite#1{\kern-20pt\hfl{#1}{}{12mm}\kern-20pt}
$$\diagram{
\azeroo{X_L}&\ldroite{f_*}&
\azeroo{X}&\ldroite{f^*}&
\azeroo{X_L}\phantom{,}\cr
\versbas{\Phi}&&\versbas{\Phi}&&\versbas{\Phi}\cr
\extun_L(\Pbar,\Lbaretoile)&\droite{f_*}&
\extun_K(\Pbar,\Lbaretoile)&\droite{f^*}&
\extun_L(\Pbar,\Lbaretoile),\cr
}\numeroter \newcount\diagext \global\diagext=\formuleno
$$}\noindent
in which the vertical arrows are the characteristic homomorphisms
of \citer\ctssequel().  

\remark{4.1} The lower $f^*$ in the diagram \formule\diagext\
consists in considering a short exact sequence of $\Gamma_K$-modules
as a short exact sequence of $\Gamma_L$-modules. The lower $f_*$ in
\formule\diagext\ associates, to the class in
$\extun_L(\Pbar,\Lbaretoile)$ of a short exact sequence
$\zero\fleche\Lbaretoile\fleche E\fleche\Pbar\fleche\zero$ of
$\Gamma_L$-modules, the class in $\extun_K(\Pbar,\Lbaretoile)$ of the
last row 
of the following
diagram of $\Gamma_K$-modules :
{\def\droite#1{\kern-5pt\hfl{#1}{}{8mm}\kern-5pt}
\def\ldroite#1{\kern-20pt\hfl{#1}{}{12mm}\kern-20pt}
$$
\diagram{
\zero&\droite{}&
\Ind\Lbaretoile&\droite{}&
\Ind E&\droite{}&
\Ind\Pbar&\droite{}&
\zero\cr
&&\pafl{}{}{5mm}&&\vershaut{}&&\vershaut{}\cr
\zero&\droite{}& 
\Ind\Lbaretoile&\ldroite{}& 
E'&\droite{}&
\Pbar&\droite{}&
\zero\cr
&&\versbas{}&&\versbas{}&&\pafl{}{}{5mm}\cr
\zero&\droite{}&
\Lbaretoile&\ldroite{}&
E''&\ldroite{}&\Pbar&\droite{}&\zero\cr}
$$}\noindent 
(defining $E'$ and $E''$), where $\Ind E$
stands for the $\Gamma_K$-module induced from the $\Gamma_L$-module
$E$, etc. ; the maps $\Ind\Lbaretoile\fleche\Lbaretoile$
(defining the last row) and $\Pbar\fleche \Ind\Pbar$ (defining the
middle row) are the natural ones \citer\serrecg(p.~13), arising from the
fact that the $\Gamma_L$-modules $\Lbaretoile$, $\Pbar$ are
restrictions of the $\Gamma_K$-modules $\Lbaretoile$, $\Pbar$.

\remark{4.2} The projection formulas
\formule\projform\ express the commutativity of certain diagrams,
namely \formule\diagcomm\ in the special case when $K$ is a finite
extension of $\qp$. We have seen that in this special case, if $X$ is
a $K$-surface potentially birational to $\P_2$, the two squares of
\formule\diagcomm\ are formally equivalent, respectively, to the two squares
\formule\diagext.  I do not know what the relationship between
\formule\projform\ and \formule\diagext\ is in the general case,
when $K$ is any perfect field and $X$ is an arbitrary $K$-variety.

\remark{4.3}  One obtains a reformulation of prop.~1.1 by
replacing, in the hypotheses, the surjectivity of $f^*$
\formule\resdef\  (resp.~$f_*$ \formule\coresdef) on reduced Brauer
groups by the injectivity of~$f_*$ (resp.~$f^*$) \formule\diagext\
on\/ $\Ext$ groups (cf.~th.~6.1).

%
%
%

\section{The general case}

As we have seen, the proof (\S\thinspace2) of prop.~1.1 hinges on the
injectivity statements of \citer\hilbertxc() and on the commutativity
--- a consequence of the projection formulas \formule\projform\ ---
of the two squares in \formule\diagcomm.  The commutativity of the
analogous squares \formule\diagext\ is a general phenomenon, as we
shall show now.  We thus obtain not only an alternate proof of
prop.~1.1 but also its extension --- along with the extensions of
prop.~3.1, cor.~3.4 and remark~3.5 which follow from it --- to the
number field case (\S\thinspace6).

\th THEOREM 5.1
\enonce
Let\/ $L\extensionde K$ be a finite extension, with\/ $K$ a perfect field and
let\/ $X$ be a smooth, proper, absolutely connected\/ $K$-variety.  Then
the two squares in diagram \formule\diagext\ are commutative.
\endth
\proof The definition of the characteristic homomorphism --- recalled
below --- makes the commutativity of the second square of
\formule\diagext\ obvious ; let us show it for the first one. 

Let $\xi$ be a $0$-cycle of degree~$0$ on $X_L$.  We then have the
$0$-cycles $\eta=f_*(\xi)$ on $X$ and $\xibar$, $\etabar$ on $\Xbar$ ;
in fact $\etabar=n\xibar$, where $n=[L:K]$.  We have to show that
$f_*(\Phi(\xi))=\Phi(\eta)$.

Let $\Abar$ be the semilocal ring of $\Xbar$ along the $0$-cycle
$\xibar$ and 
put $\Ubar=U\times_K\Lbar$, where $U\subset X$ is the
open subvariety complementary to the support of $\eta$.  
Recall \citer\ctssequel() that $\Phi(\xi)\in\extun_L(\Pbar,\Lbaretoile)$ is
the class of the extension 
of $\Gamma_L$-modules deduced from the short exact sequence of
$\Gamma_K$-modules 
$$
\one\fleche\Abaretoile\!/\,\Lbaretoile\fleche
\hbox{\ss Z}^1(\Ubar)\fleche
\Pbar\fleche\zero
\numeroter \newcount\theext \global\theext=\formuleno
$$ 
via the $\Gamma_L$-homomorphism of evaluation
$\Abaretoile\!/\,\Lbaretoile\fleche\Lbaretoile$ at $\xibar$ ;
similarly, $\Phi(\eta)\in\extun_K(\Pbar,\Lbaretoile)$ is deduced
from \formule\theext\ via the $\Gamma_K$-homomorphism of evaluation
$\Abaretoile\!/\,\Lbaretoile\fleche\Lbaretoile$ at $\etabar$.  For
$f_*(\Phi(\xi))$, see the explicit description of the map $f_*$
given in remark~4.1.  To show that $f_*(\Phi(\xi))=\Phi(\eta)$, we use
the following lemma.

\th LEMMA 5.2
\enonce
Let $P$, $Q$, $R$ be $\Gamma_K$-modules and let $h:Q\fleche R$ be a
$\Gamma_L$-homomorphism such that $nh:Q\fleche R$ is
$\Gamma_K$-equivariant.  Then the following square is commutative :
$$
\diagram{
\extun_K(P,Q)&\droite{nh}&
\extun_K(P,R)\phantom{.}\cr
\versbas{f^*}&&\vershaut{f_*}\cr
\extun_L(P,Q)&\droite{h}&
\extun_L(P,R).\cr
}$$
\endth


\section{The global application}

With the general result (th.~5.1) at hand, the statements and proofs
of \S\S\thinspace2,3 can be extended to any situation where the
characteristic homomorphisms are injective.  We shall content
ourselves with the following

\th THEOREM 6.1 
\enonce
Assume that\/ $X$ is a\/ $K$-surface\/ potentially birational to
$\P_2$, where\/ $K\extensionde \Q$ is either a finite extension or a
finitely generated extension and, in the latter case, that\/ $X$ has
a\/ $K$-rational\/ $0$-cycle of degree~$1$.  Let\/ $L\extensionde K$
be a finite extension.

a) If the corestriction~$f_*$ \formule\diagext\ on\/ $\Ext$ groups
is injective (resp.~$0$), then the corestriction~$f_*$
\formule\coresdef\ on reduced Chow groups is injective (resp.~$0$).

b) If the restriction~$f^*$ \formule\diagext\ on\/ $\Ext$ groups is
injective (resp.~$0$), then the restriction~$f^*$
\formule\resdef\ on reduced Chow groups is injective (resp.~$0$).

%
%
\endth
\proof The ingredients of the proof are the same as in the local case
(\S\S\thinspace2,3), up to replacing the diagram \formule\diagcomm\
by the diagram \formule\diagext.

It was shown in \citer\hilbertxc() that the
vertical arrows in \formule\diagext\ --- the characteristic
homomorphisms of \citer\ctssequel() --- are injective under our
hypotheses ; the result follows therefrom and the commutativity of
\formule\diagext\ (th.~5.1).



\unvbox\bibbox 
\vskip2cm
{\obeylines\parskip=0pt\parindent=0pt
Chandan Singh Dalawat
Harish-Chandra Research Institute
Chhatnag Road, Jhunsi
Allahabad 211\thinspace019, India.
\line{\hfill}
\tt dalawat@mri.ernet.in}
\bye